\documentclass{smfart}
\usepackage{graphicx}
\usepackage{amssymb}
\usepackage{epstopdf}

\usepackage{amsmath}
\usepackage{amsthm}
\usepackage{pdfsync}
\usepackage{hyperref}
\usepackage{color}

\newcommand{\bb}[1]{\mathbb{#1}}

\newtheorem{theorem}{Theorem}

\def\zz{\mathbb{Z}}
\def\rr{\mathbb{R}}

\def\eps{\varepsilon}

\def\bu{{\bf{u}}}

\title[Dispersion for Schr\"odinger  equations]{ The Dispersion property for Schr\"odinger  equations}
\author{Liviu I. Ignat}

\address{L. I. Ignat
\hfill\break\indent Institute of Mathematics ``Simion Stoilow'' of the Romanian Academy\\
\hfill\break\indent  21 Calea Grivitei Street \\010702 Bucharest \\ Romania 
\hfill\break\indent \and
\hfill\break\indent Faculty of Mathematics and Computer Science, University of Bucharest \\
\hfill\break\indent 14 Academiei Street, 010014 Bucharest, Romania.
}
 
 \email{{\tt
liviu.ignat@gmail.com}\hfill\break\indent  {\it Web page: }{\tt
http://www.imar.ro/\~\,lignat}}

\begin{document}


\begin{abstract}
In this paper we analyze the dispersion property of some models   involving Schr\"odinger equations. First we  focus on the discrete case and then we present some results on graphs.
\end{abstract}

\maketitle

\section{Introduction}

Let us first recall some classical properties of the solutions of the Schr\"odinger equation. The solution of the equation
\begin{equation}\label{sch-liniar}
\left\{
\begin{array}{ll}
iu_t+u_{xx}=0,& x\in \rr,t\in \rr,\\
u(0,x)=\varphi(x), & x\in \rr,
\end{array}
\right.
\end{equation}
can be obtained by using the Fourier transform as follows
\[
u(t)=S(t)\varphi = (e^{-4\pi^2 it|\xi|^2}\widehat{\varphi})^\vee=\frac{e^{i|x|^2/4t}}{(4\pi i t)^{1/2}}\ast \varphi.
\] 
There are two properties that follow from the above representation. The first one is the conservation of the $L^2(\rr)$ norm:
\begin{equation}\label{l.2}
\|S(t)\varphi\|_{L^2(\rr)}=\|\varphi\|_{L^2(\rr)}.
\end{equation}
The second one is the so-called  \textit{dispersive} property, which shows that the solutions of system \eqref{sch-liniar} decay as time increases:
\begin{equation}\label{decay.infty}
\|S(t)\varphi\|_{L^\infty(\rr)}\leq \frac{1}{(4\pi |t|)^{1/2}}\|\varphi\|_{L^1(\rr)}.
\end{equation}
These simple properties can be used  to obtain more refined estimates for the linear semigroup. There are  properties of gain on integrability/regularity with respect to the initial data: the Strichartz estimates and the so-called local smoothing property. We state them here in their simplest form. Similar estimates can be written for the inhomogenous problem by using the $TT^*$ argument and eventually Christ-Kiselev's argument \cite{0974.47025}. For any initial data in $L^2(\rr)$ Strichartz's estimates state that the solutions of \eqref{sch-liniar} satisfy 
\begin{equation}\label{strichartz}
\|S(t)\varphi\|_{L^q_t(\rr,L^r_x(\rr))}\lesssim \|\varphi\|_{L^2(\rr)}.
\end{equation}
A scaling argument forces the following admissibility condition for the pairs $(q,r)$:
\[
\frac 2q+\frac 1r=\frac 12.
\]
The local smoothing property means that the linear semigroup  gains locally a one-half space derivative with respect to the initial data:
\begin{equation}\label{gail.derivative}
\sup_{x\in \rr} \int _{\rr} |(-\Delta)^{1/4} (S(t)\varphi) |^2dt \lesssim \|\varphi\|_{L^2(\rr)}.
\end{equation}

The above estimates play a crucial role in proving the well-posedness of the nonlinear Schr\"odinger equation. We recall here only one result in this direction \cite{0638.35021} and the interested reader should consult the classical references \cite{1055.35003,MR2233925,MR2492151,MR1696311}.   For any  initial data in $L^2(\rr)$ and $p\in (0,4)$ there exists a unique solution $u\in C(\rr,L^2(\rr))\cap L^q_{loc}(\rr,L^r(\rr))$ of the equation
\begin{equation}\label{sch-nonlinear}
\left\{
\begin{array}{ll}
iu_t+u_{xx}=|u|^pu,& x\in \rr,t\in \rr,\\
u(0,x)=\varphi(x), & x\in \rr.
\end{array}
\right.
\end{equation}

In this paper we analyze the  dispersion property \eqref{decay.infty} in the context of  discrete equations and of the equations on graphs. Similar estimates in the spirit of \eqref{strichartz} and \eqref{gail.derivative} can also be obtained. 

\section{Discrete equations}
Let us consider the discrete Schr\"odinger equation 
\begin{equation}\label{sch-discrete}
\left\{
\begin{array}{ll}
iu_t+\Delta_d u=0,& j\in \zz,t\in \rr,\\
u(0,j)=\varphi(j), & j\in \zz,
\end{array}
\right.
\end{equation}
where, for any  function $u:\zz\rightarrow \rr$, the discrete Laplacian $\Delta _d$ acts as follows:
\[
(\Delta_d u)(j)=u(j+1)-2u(j)+u(j-1).
\]
In this paper we will not discuss  the possible numerical approximation of the NSE as  in \cite{MR2485456, ignat-zuazua-jmpa, MR3119633, MR2755707} where there is a parameter $h$, the mesh size, that it is going to zero.

In the case of the discrete equation \eqref{sch-discrete} the solutions decay as follows:
\begin{theorem}\label{stefanov}(\cite{MR2135236, MR2150357})
For any $\varphi\in l^1(\zz)$ the solution of system \eqref{sch-discrete} satisfies
\begin{equation}\label{decay.discrete}
\|u(t)\|_{l^\infty(\zz)}\lesssim (|t|+1)^{-1/3}\|\varphi \|_{L^1(\zz)}, \quad \forall \ t\in \rr.
\end{equation}
\end{theorem}
Similar to the case of the real line we can write the solution of system \eqref{sch-discrete} as a convolution:
\[
u(t,j)=\sum _{k\in \zz}K_t(j-k)\varphi(k),
\]
where
\[
K_t(j)=\int _{-\pi}^{\pi} e^{-4it \sin ^2(\xi/2)}e^{ij\xi}d\xi.
\]
Using a classical result for oscillatory integrals, the Van der Corput Lemma \cite[Ch. 1]{MR2492151}, we have that 
\[
|K_t(j)|\lesssim (1+|t|)^{-1/3}, \quad \forall j\in\zz .
\]
This is a consequence of the fact that the second and the third derivatives of the phase
 function $\psi(\xi)=4\sin ^2(\xi/2)$ do not vanish simultaneously at the same point: 
$|\psi'' |+|\psi '''|\geq C>0$. Observe that $\|K_t\|_{L^2(\zz)}\leq C$ for some positive constant, so, for any $p\in [2,\infty]$, 
\[
\|K_t\|_{l^p(\zz)}\lesssim (|t|+1)^{-1/3}, \quad \forall \ t\in \rr.
\]
This implies a  decay property for $u(t)$ in the $l^p(\zz)$-norm similar to the one in \eqref{decay.discrete}
$
\|u(t)\|_{L^p(\zz)}\lesssim (|t|+1)^{-(p-2)/3p}\|\varphi \|_{L^1(\zz)}.
$ 
However, in \cite{MR2682115} the authors show that this decay property can be improved as follows
\[
\|u(t)\|_{L^p(\zz)}\lesssim (|t|+1)^{-\alpha_p}\|\varphi \|_{L^1(\zz)},
\] 
where
\[
\alpha_p=\left\{
\begin{array}{rr}
\frac{p-2}{2p}, & p\in [2,4),\\[10pt]
\frac{p-1}{3p}, & p\in (4,\infty].
\end{array}
\right.
\]
Discrete equations on  positive integers can also be considered
\begin{equation}\label{dlse}
\left\{
\begin{array}{ll}
  i u_t(t,j) + (\Delta _du)(t,j) = 0,&j\geq 1, t\neq 0,\\
   u(0,j) = \varphi(j), &j\geq 1.
\end{array}
\right.
\end{equation}
One may also impose 
Dirichlet or Neumann boundary conditions at $j=1$: $u(t,0)=0$ or $u(t,0) = u(t,1)$, respectively. The solutions of these systems have similar decay.

In a joint paper with D. Stan \cite{ignat-diana} we considered the case of two discrete Schr\"odinger equations coupled at $j=0$:
 \begin{equation}\label{coupled-diana}
\left\{
\begin{array}{ll}
  i u_t(j) + b_1^{-2}(\Delta_d u)(j) = 0& j\leq -1, \\
i v_t(j) + b_2^{-2}(\Delta_d v)(j) = 0& j\geq 1, \\
  u(t,0) = v(t,0) , &t>0, \\
  b_1^{-2}(u(t,-1)-u(t,0))=b_2^{-2}(v(t,0)-v(t,1)),&t>0,\\
  u(0,j) = \varphi(j),& j \leq -1,\\
    v(0,j) = \varphi(j),& j \geq 1.
\end{array}
\right.
\end{equation}

In the particular case  when $b_1=b_2$ the solution of system \eqref{coupled-diana} can be written in terms of the solutions of some Dirichlet and Neumann problems. Thus, it decays  with the same power $(|t|+1)^{-1/3}$. For $b_1\neq b_2$ there are some difficulties.  
In matrix formulation system \eqref{coupled-diana} can be written as $iU_t+AU=0$ where  $U=(v(-j),u(j))_{j\neq 0}$  and 
 \[
 A=\left(\begin{array}{cccccccc}
... & ... & ... & 0 & 0 & 0 & 0 & 0 \\
0 & b_1^{-2} & -2b_1^{-2} & b_1^{-2} & 0 & 0 & 0 & 0 \\
0 & 0 & b_1^{-2} & - b_1^{-2}-\frac{1}{b_1^2+b_2^2}& \frac{1}{b_1^2+b_2^2}& 0 & 0 & 0 \\
0 & 0 & 0 & \frac{1}{b_1^2+b_2^2}& -\frac{1}{b_1^2+b_2^2}-b_2^{-2} & b_2^{-2} & 0 & 0 \\
0 & 0 & 0 & 0 & b_2^{-2} & -2b_2^{-2} & b_2^{-2} & 0 \\
0 & 0 & 0 & 0 & 0 & ... & ... & ...
\end{array}\right).
\]
 Since the matrix is not diagonal we cannot use  the Fourier transform. We use the explicit form of the resolvent $R(\lambda)=(A-\lambda I)^{-1}$ and  spectral calculus  to write the following limiting absorption principle for $\varphi\in L^1(\rr)$
 \[
e^{itA}\varphi =\frac 1{2i\pi} \int _{\sigma (A)} e^{it\omega} (R^+(\omega)-R^-(\omega))\varphi d\omega,
\] 
 where $R^{\pm}= \lim _{\eps \downarrow 0}R(\omega \pm i \eps)$.
Using this representation we obtain that the solutions of system \eqref{coupled-diana} also decay as $(|t|+1)^{-1/3}$. To prove the  decay property it is enough to show that for any $a\in (0,1]$  there is a positive constant $C(a)$ such that 
 the following 
\begin{equation}\label{int.osc}
\Big| \int _0^\pi e^{it (2\cos \theta +2z\arcsin (a\sin \frac \theta 2))}e^{ity\theta}\sin \theta d\theta \Big|\leq C(a)(|t|+1)^{-1/3}
\end{equation}
holds for any real numbers $y$, $z$ and $t$. The main difficulty comes from the fact that we want an estimate that should be uniform with respect to the parameters $y$ and $z$. To obtain that we use previous improvements of the Van der Corput lemma given in \cite{0738.35022} and new ones.
Full details are given in \cite{ignat-diana}.
 
 Many questions can arise regarding the problems discussed above. Using matrices, the discrete equations can be rewritten as follows:
\begin{equation}\label{sys.A}
iU_t+AU=0.
\end{equation}
We can assume that matrix $A$ has a few non-identical vanishing diagonals. One of the  questions that arises here is the following: what  properties of matrix $A$  guarantee the decay of the solutions of system \eqref{sys.A}? We admit that also other types of decay may appear 
$t^{-1/4}$, $t^{-\alpha}\log(t)$, etc... In the case when $A$ is a diagonal matrix we can use Fourier analysis to reduce the question to the study of the roots of some trigonometric polynomials. It will  be  interesting to see how resolvent estimates similar to those in 
\cite{MR2227135} used in the continuous case 
can be obtained and used in the matrix case.  The analysis of system \eqref{coupled-diana} has been mainly done in the context of the $l^1(\zz)-l^\infty(\zz)$ decay of the solutions. It is still possible to refine the analysis to obtain better estimates as those in \cite{MR2682115}.
 
We have mainly analyzed the particular case of two discrete equations  coupled at one point, but more general coupling conditions  can  be used or more coupled equations can be considered. We recall here \cite{gavrus} where some discrete models on trees have been analyzed.

 \section{Equation on graphs}\label{graphs}
 
 In this section we present some results concerning the Schr\"odinger equation on graphs. The results obtained until now concern trees, i.e. graphs without cycles, having the last generation of edges formed from infinite strips. For a graph $\Gamma=(V,E)$ we consider the following equation
 \begin{equation}\label{eq.tree}
\left\{
\begin{array}{ll}
i\bu _t(t,x)+\Delta_\Gamma \bu(t,x)=0,& x\in \Gamma, t\neq 0 ,\\[10pt]
\bu(0,x)=\bu_0(x),&   x\in \Gamma.
\end{array}
\right.
\end{equation}
For the precise definition of operator $\Delta_\Gamma$ we refer to \cite{MR2684310}. Denoting $\bu=(u_e)_{e\in E}$
system \eqref{eq.tree} means that on each edge $e\in E$ of the graph we have a Schr\"odinger  equation for $u_e$ and at any  vertex $v\in V$   we couple the equations on the edges that enter into $v$
 by assuming for example,  Kirchhoff's law: continuity and sum of the normal derivatives to be zero. Other coupling conditions can be imposed:  $\delta$ or $\delta'$ coupling as in \cite{MR2804557} or those given in \cite{MR1671833}. 

The first result regarding the dispersion property on trees has been obtained in \cite{MR2684310} where the case of regular trees has been considered. The regular trees are special trees  where all the edges in the same generation have the same length and all the vertices in the same generation have the same number of children. In this case we can make averages of the  solutions defined in the same generation and reduce the analysis of dispersion on the tree to the case of the Schr\"odinger equation on the real line with a piecewise constant coefficient $\sigma$, $iu_t+(\sigma u_x)_x=0$. When the coefficient $\sigma$ is given by a finite number of piecewise constant functions we have the following result.
\begin{theorem}[\cite{MR2049025}] Consider a partition of the real axis 
$-\infty=x_0<x_1<\dots<x_{n+1}=\infty$
and a step function 
$\sigma(x)=\sigma_i \ \text{for}\ x\in (x_i,x_{i+1}),$
where $\sigma_i$ are positive numbers.
The solution $u$ of the Schr\"odinger equation 
\[
\left\{
\begin{array}{ll}
iu_t(t,x)+(\sigma(x)u_x)_x(t,x)=0,& \text{for}\ x\in \rr,t\neq 0,\\[10pt]
u(0,x)=\varphi(x), &x\in \rr,
\end{array}
\right.
\]
satisfies the dispersion inequality
\[
\|u(t,\cdot)\|_{L^\infty(\rr)}\leq C|t|^{-1/2}\|\varphi\|_{L^1(\rr)}, \quad \ t\neq 0.
\]
\end{theorem}
Using this result and an inductive argument on the  generations of the tree, in \cite{MR2684310} the dispersion property for the solutions of system [13] has been proved.

The case of the star-shaped tree has  been considered  in \cite{MR2804557} not only with the classical Kirchhoff's coupling, but also with the $\delta$ and $\delta'$ coupling. The advantage of the star-shaped tree is  that explicit formulas can be easily obtained for the resolvent $(\Delta_\Gamma-\lambda I)^{-1}$ and then for the solutions of system \eqref{eq.tree}. Classical results for oscillatory integrals allow us to obtain  the dispersion property more easily.

The case of general trees with Kirchhoff's coupling has been considered in \cite{ignat-banica}. In contrast to the case of the star-shaped tree, the resolvent is written in an implicit way by using results for almost periodic functions. The same problem with $\delta$-coupling has been treated in \cite{MR3254348}. We emphasize that  in the particular  case of a tree with all internal vertices having degree two  this type of coupling corresponds to  the Schr\"odinger equation on the real line where the Laplacian has been perturbed by a sum of Dirac deltas: 
\[
H_{\alpha}=-\Delta + \sum _{j=1}^p \alpha_j \delta (x-x_j).
\]
Under some technical assumptions  on the strengths $\{\alpha_j\}_{j=1}^p$ and on the lengths of the intervals $\{x_{j+1}-x_j\}_{j=1}^{p-1}$ that exclude the existence of resonances (\cite{MR2608273,MR2791127,MR3254348}) it has been proved that
\begin{equation}\label{dispersion-0}
\| e^{-itH_\alpha}Pu_0\|_{L^\infty(\rr)}\leq \frac{C}{\sqrt {|t|}}
\| u_0\|_{L^1(\rr)},\,\,\, \forall t\neq 0,
\end{equation}
where $P$ is the projection on the continuous part of the spectrum of $H_\alpha$. A similar analysis can be performed on a tree 
with delta coupling conditions  at the vertices.
The interested reader cans check all the details in \cite{MR3254348}.
 
 The analysis offered so far  concerns Kirchhoff's coupling or the $\delta$-coupling. The $\delta'$-coupling has been only considered  in the case of a star-shaped tree \cite{MR2804557}. In the case of metric graphs, there are various types of coupling conditions that introduce a self-adjoint version of the Laplace operator $\Delta(A,B)$. In \cite{MR1671833} it was proved that for the coupling
\[
A(v){ \bu } (v)+B(v){\bu }'(v)=0,
\]
the  following conditions are necessary and sufficient:
 $(A(v), B(v))$ has maximal rank, i.e.  $d(v)$, and $A(v)B(v)^T=B(v)A(v)^T$ for all vertices $v$ of the graph $\Gamma$.
As far as we know up to now, the dispersion property has not been proved under general coupling conditions  not even on star shaped trees. 

One of the problems that emerges from the previous analysis is how to obtain  some kind of dispersion when the graphs have  cycles.  We recall that in the case of a compact  manifold without boundary \cite{MR2058384} is was proved that the dispersion property holds for small intervals of time depending on the range of the  eigenfunctions taken by  the initial data. In the case of  one-dimensional torus $\mathbb{T}^1$ this was already  known \cite[Th. 5.3]{0738.35022}: 
\[
\Big\|\sum _{|k|\leq N} a_k e^{i(t k^2+kx)} \Big\| _{L^\infty(\mathbb{T}^1)}\lesssim |t|^{-1/2} \Big\|\sum _{k} a_k e^{ikx} \Big\|_{L^1(\mathbb{T}^1)}, \quad \forall |t|\leq N^{-1}.
\]
The existence of   an argument based on oscillatory integrals that allows us to obtain the dispersion for $\mathbb{T}^1$  on small time intervals gives us  hope that in the case of graphs with cycles similar results could be obtained. In fact, for very particular structures as those in \cite{MR1207868} the analysis performed in \cite{grecu} shows that the dispersion holds for small intervals of time. The full understanding of the dispersion phenomena when a periodic structure is combined with  infinite lines,  i.e. a graph with few infinite external edges that contains a cycle, remains to be investigated.  In this context we recall the case of a cylinder that was considered in \cite{0976.35085}.

Another model that can be considered on a graph structure is Dirac's equation:
$
iu_t=\mathcal{H}u
$
 where 
\[
\mathcal{H}=\left[\begin{array}{cc}-i \partial_x & -1 \\
-1 & i\partial_x
\end{array}\right].
\]
 As far as we know the dispersion property for this equation on graphs/trees has not been considered previously.

\medskip
 {\bf
Acknowledgements.}

L. I.  Ignat was partially supported by Grant PN-II-ID-PCE-2011-3-0075 of the Romanian National Authority for Scientific Research, CNCS -- UEFISCDI
and by LEA Franco-Roumain Math-Mode. The author thanks the organizers of the workshop 
"Pde's, Dispersion, Scattering theory and Control theory", Monastir, 2013,  for their kind invitation.  
The author also thanks Scoala Normala Superioara Bucuresti and his collaborators  V. Banica, C. Gavrus, A. Grecu and D. Stan.


\begin{thebibliography}{10}

\bibitem{MR2804557}
R.~Adami, C.~Cacciapuoti, D.~Finco, and D.~Noja.
\newblock Fast solitons on star graphs.
\newblock {\em Rev. Math. Phys.}, 23(4):409--451, 2011.

\bibitem{MR3119633}
C.~ Audiard.
\newblock Dispersive schemes for the critical {K}orteweg-de {V}ries equation.
\newblock {\em Math. Models Methods Appl. Sci.}, 23(14):2603--2646, 2013.

\bibitem{MR2049025}
V.~Banica.
\newblock Dispersion and {S}trichartz inequalities for {S}chr\"odinger
  equations with singular coefficients.
\newblock {\em SIAM J. Math. Anal.}, 35(4):868--883 (electronic), 2003.

\bibitem{MR3254348}
V.~Banica and L.I.~Ignat.
\newblock Dispersion for the {S}chr\"odinger equation on the line with multiple
  {D}irac delta potentials and on delta trees.
\newblock {\em Anal. PDE}, 7(4):903--927, 2014.

\bibitem{ignat-banica}
Valeria Banica and L.I.~Ignat.
\newblock Dispersion for the schr\"odinger equation on networks.
\newblock {\em J. Math. Phys.}, 52:083703, 2011.

\bibitem{MR2058384}
N.~Burq, P.~G{\'e}rard, and N.~Tzvetkov.
\newblock Strichartz inequalities and the nonlinear {S}chr\"odinger equation on
  compact manifolds.
\newblock {\em Amer. J. Math.}, 126(3):569--605, 2004.

\bibitem{MR2227135}
N.~Burq and F.~Planchon.
\newblock Smoothing and dispersive estimates for 1{D} {S}chr{\"o}dinger
  equations with {BV} coefficients and applications.
\newblock {\em J. Funct. Anal.}, 236(1):265--298, 2006.

\bibitem{1055.35003}
T.~Cazenave.
\newblock {\em {Semilinear Schr\"{o}dinger equations.}}
\newblock {Courant Lecture Notes in Mathematics 10. Providence, RI: American
  Mathematical Society (AMS); New York, NY: Courant Institute of Mathematical
  Sciences. xiii }, 2003.

\bibitem{0974.47025}
M.~Christ and A.~Kiselev.
\newblock {Maximal functions associated to filtrations.}
\newblock {\em J. Funct. Anal.}, 179(2):409--425, 2001.

\bibitem{MR2791127}
V.~Duch{\^e}ne, J.~L. Marzuola, and M.I. Weinstein.
\newblock Wave operator bounds for one-dimensional {S}chr\"odinger operators
  with singular potentials and applications.
\newblock {\em J. Math. Phys.}, 52(1):013505, 17, 2011.

\bibitem{MR1207868}
B.~Gaveau, M.~Okada, and T.~Okada.
\newblock Explicit heat kernels on graphs and spectral analysis.
\newblock In {\em Several complex variables ({S}tockholm, 1987/1988)},
  volume~38 of {\em Math. Notes}, pages 364--388. Princeton Univ. Press,
  Princeton, NJ, 1993.

\bibitem{gavrus}
C.~Gavrus.
\newblock Dispersion property for discrete {S}chr\"odinger equations on
  networks.
\newblock Master's thesis, Scoala Normala Superioara Bucuresti, July 2012.

\bibitem{grecu}
A.~Grecu.
\newblock Dispersion property for {S}chr\"odinger equations on graphs with
  cycles.
\newblock Master's thesis, Univ. of Bucharest, in preparation.

\bibitem{MR2684310}
L.I.~Ignat.
\newblock Strichartz estimates for the {S}chr{\"o}dinger equation on a tree and
  applications.
\newblock {\em SIAM J. Math. Anal.}, 42(5):2041--2057, 2010.

\bibitem{ignat-diana}
L.I.~Ignat and D.~Stan.
\newblock Dispersive properties for discrete schr{\"o}dinger equations.
\newblock {\em Journal of Fourier Analysis and Applications}, 17(5):1035--1065,
  2011.

\bibitem{MR2135236}
L.I.~Ignat and E.~Zuazua.
\newblock Dispersive properties of a viscous numerical scheme for the
  {S}chr{\"o}dinger equation.
\newblock {\em C. R. Math. Acad. Sci. Paris}, 340(7):529--534, 2005.

\bibitem{MR2485456}
L.I.~Ignat and E.~Zuazua.
\newblock Numerical dispersive schemes for the nonlinear {S}chr{\"o}dinger
  equation.
\newblock {\em SIAM J. Numer. Anal.}, 47(2):1366--1390, 2009.

\bibitem{ignat-zuazua-jmpa}
L.I.~Ignat and E.~Zuazua.
\newblock Convergence rates for dispersive approximation schemes to nonlinear
  schr\"odinger equations.
\newblock {\em J. Math. Pures Appl.}, to appear.

\bibitem{0738.35022}
C.E. Kenig, G.~Ponce, and L.~Vega.
\newblock {Oscillatory integrals and regularity of dispersive equations.}
\newblock {\em Indiana Univ. Math. J.}, 40(1):33--69, 1991.

\bibitem{MR1671833}
V.~Kostrykin and R.~Schrader.
\newblock Kirchhoff's rule for quantum wires.
\newblock {\em J. Phys. A}, 32(4):595--630, 1999.

\bibitem{MR2608273}
H.~Kova{\v{r}}{\'{\i}}k and A.~Sacchetti.
\newblock A nonlinear {S}chr{\"o}dinger equation with two symmetric point
  interactions in one dimension.
\newblock {\em J. Phys. A}, 43(15):155205, 16, 2010.

\bibitem{MR2492151}
F.~Linares and G.~Ponce.
\newblock {\em Introduction to nonlinear dispersive equations}.
\newblock Universitext. Springer, New York, 2009.

\bibitem{MR2755707}
A.~Marica and E.~Zuazua.
\newblock High frequency wave packets for the {S}chr\"odinger equation and its
  numerical approximations.
\newblock {\em C. R. Math. Acad. Sci. Paris}, 349(1-2):105--110, 2011.

\bibitem{MR2682115}
A.~Mielke and C.~Patz.
\newblock Dispersive stability of infinite-dimensional {H}amiltonian systems on
  lattices.
\newblock {\em Appl. Anal.}, 89(9):1493--1512, 2010.

\bibitem{MR2150357}
A.~Stefanov and P.G. Kevrekidis.
\newblock Asymptotic behaviour of small solutions for the discrete nonlinear
  {S}chr\"{o}dinger and {K}lein-{G}ordon equations.
\newblock {\em Nonlinearity}, 18(4):1841--1857, 2005.

\bibitem{MR1696311}
C.~Sulem and P.L. Sulem.
\newblock {\em The nonlinear {S}chr\"odinger equation}, volume 139 of {\em
  Applied Mathematical Sciences}.
\newblock Springer-Verlag, New York, 1999.
\newblock Self-focusing and wave collapse.

\bibitem{0976.35085}
H.~Takaoka and N.~Tzvetkov.
\newblock {On 2D nonlinear Schr{\"o}dinger equations with data on $\bb{R}\times
  \bb{T}$.}
\newblock {\em J. Funct. Anal.}, 182(2):427--442, 2001.

\bibitem{MR2233925}
T.~Tao.
\newblock {\em Nonlinear dispersive equations}, volume 106 of {\em CBMS
  Regional Conference Series in Mathematics}.
\newblock Published for the Conference Board of the Mathematical Sciences,
  Washington, DC, 2006.
\newblock Local and global analysis.

\bibitem{0638.35021}
Y.~Tsutsumi.
\newblock {$L\sp 2$-solutions for nonlinear Schr\"{o}dinger equations and
  nonlinear groups.}
\newblock {\em Funkc. Ekvacioj, Ser. Int.}, 30:115--125, 1987.

\end{thebibliography}
%

\end{document}